\pdfoutput=1
\documentclass{amsart}

\DeclareMathSymbol{\twoheadrightarrow}  {\mathrel}{AMSa}{"10}
        \def\GG{{\mathcal G}}

\def\Q{{\mathbb Q}}
\def\Z{{\mathbb Z}}
\def\C{{\mathbb C}}
\def\CC{{\mathfrak C}}
\def\RR{{\mathbb R}}
\def\F{{\mathbb F}}
\def\P{{\mathbb P}}

\def\f{{\tilde F}}
                     \def\f0{{\mathfrak f}}

             \def\K{\mathrm{K}}

\def\A8{{\mathbf A}_8}
\def\Alt{\mathrm{Alt}}
\def\RR{{\mathfrak R}}
\def\Perm{\mathrm{Perm}}
\def\Gal{\mathrm{Gal}}
\def\Pic{\mathrm{Pic}}

\def\End{\mathrm{End}}
\def\Aut{\mathrm{Aut}}
\def\div{\mathrm{div}}

              \def\L{{\mathcal L}}

\def\ST{{\mathbf S}}

\def\fchar{\mathrm{char}}
\def\GL{\mathrm{GL}}

\def\M{\mathrm{M}}
\def\A{\mathbf{A}}

\def\B{{\mathcal B}}

\def\OO{\mathrm{O}}
            \def\SO{\mathrm{SO}}
\def\dim{\mathrm{dim}}
\def\Oc{{\mathcal O}}

       \def\ZZ{{\mathcal Z}}
       \def\PGL{\mathrm{PGL}}

\newtheorem{thm}{Theorem}[section]
\newtheorem{lem}[thm]{Lemma}

\newtheorem{prop}[thm]{Proposition}
\theoremstyle{definition}
\newtheorem{defn}[thm]{Definition}
\newtheorem{ex}[thm]{Example}
\newtheorem{rem}[thm]{Remark}

        \newtheorem{sect}[thm]{}
\title[Del Pezzo surfaces and jacobians]
{Del Pezzo surfaces of degree $1$ and  jacobians}

\author[Yu.\ G.\ Zarhin]{Yu.\ G.\ Zarhin}
\address{Department of Mathematics, Pennsylvania State University,
University Park, PA 16802, USA}
\email{zarhin\char`\@math.psu.edu}
\thanks{Supported by SFB 701 ``Spektrale Strukturen und topologische Methoden in der Mathematik"
 (Fakult\"at f\"ur Mathematik der Universit\"at Bielefeld)}
\begin{document}

\maketitle

\section{Introduction}
Let $K$ be a field of characteristic zero, $\bar{K}$ its algebraic
closure and $\Gal(K)=\Aut(\bar{K}/K)$ its absolute Galois group.

In \cite{ZarhinMRL} the author constructed explicitly
$g$-dimensional abelian varieties (jacobians) without non-trivial
endomorphisms for every $g>1$. This construction may be described
as follows. Let $n=2g+1$ or $2g+2$. Let us choose an $n$-element
set $\RR \in \bar{K}$ that constitutes a Galois orbit over $K$ and
assume, in addition, that the Galois group of $K(\RR)$ over $K$
coincides either with the full symmetric group $\ST_n$ or the
alternating group $\A_n$. Let $f(x) \in K[x]$ be the irreducible
polynomial of degree $n$, whose set of roots coincides with $\RR$.
Let us consider the genus $g$ hyperelliptic curve $C_f:y^2=f(x)$
over $\bar{K}$ and let $J(C_f)$ be its jacobian, which is the
$g$-dimensional abelian variety. Then the ring $\End(J(C_f))$ of
all $\bar{K}$-endomorphisms of $J(C_f)$ coincides with $\Z$.

It is well-known that every genus $2$ curve is hyperelliptic.
However, there is a plenty of non-hyperelliptic genus $3$ curves:
namely, a curve of genus $3$ is  non-hyperelliptic if and only if
it is isomorphic to a smooth plane quartic. So, one may ask for a
natural construction of such quartics, whose jacobians have no
nontrivial endomorphisms. In order to do that, suppose that we are
given seven $\bar{K}$-points on the projective plane {\sl in
general position}, i.e., no three points lie on a one line and no
six on a one conic.  Assume also that $\Gal(K)$ permutes those
seven points transitively. By blowing them up, we obtain a Del
Pezzo surface of degree $2$ that is defined over $K$ (\cite[\S
3]{Iskovskikh}, \cite[Th. 1 on p. 27]{Demazure}). Suppose that the
$7$-element Galois orbit has {\sl large} Galois group, namely,
either  $\ST_7$ or the alternating group $\A_7$. It is proven in
\cite{ZarhinVE} that if we consider the anticanonical map of the
Del Pezzo surface onto the projective plane then the jacobian of
the corresponding branch curve  has no nontrivial endomorphisms
over $\bar{K}$. (Recall \cite[Ch. VII, Sect. 4]{Dolgachev} that
this curve is a smooth plane quartic.) Also, starting with an
irreducible degree $7$ polynomial with large Galois group, we
provided an explicit construction of a $7$-element Galois orbit in
general position and with the same Galois group. (If one
starts with an irreducible quartic polynomial $f(x)$ over $K$
 with Galois group $\ST_4$ and considers
 a smooth plane quartic $C_{f,3}:y^3=f(x)$ then it turns out that the endomorphism ring of its jacobian
  is $\Z[\frac{-1+\sqrt{-3}}{2}]$
 if $K$ contains $\sqrt{-3}$) \cite{ZarhinD}.)

The aim of this paper is  to prove a similar result while dealing
with eight points, Del Pezzo surfaces of degree $1$ and their
branch curves with respect to bi-anticanonical maps. (In this case
the curve involved is a (non-hyperelliptic) genus $4$ curve with
vanishing theta characteristic \cite{Demazure,Dolgachev}).
Notice that Del Pezzo surfaces of degree $1$ do depend on $8$ ``parameters" while the moduli space of genus $4$ curves has
dimension $9$. So, it is not {\sl apriori} obvious (at least, to the author) why there exists
(even over the field $\C$ of complex numbers) a degree $1$ Del Pezzo surface with  simple jacobian of its branch curve.

We
prove that the endomorphism algebra of the corresponding jacobian
(over $\bar{K}$) is either $\Q$ or a quadratic field; in
particular, the jacobian is an absolutely simple abelian fourfold.
Also, starting with an irreducible degree $8$ polynomial with
large Galois group (${\mathrm S}_8$ or  ${\mathrm A}_8$), we
provide an explicit construction of a $8$-element Galois orbit
with the same Galois group and in general position. (In the case
of {\sl eight} points, we have to check additionally that there is
no cubic that contains all the points and one of those points is
singular on the cubic \cite[Sect. 3]{Iskovskikh}, \cite[Th. 1 on
p. 27]{Demazure}).     In particular, we prove the following statement.

\begin{thm}
\label{mainex}
   Let $f(t)\in K[t]$ be an irreducible  degree $8$ polynomial, whose Galois
group is either $\ST_8$ or $\A_8$. Let $\RR\subset K_a$ be the set
of roots of $f$. Then:
\begin{itemize}
\item[(i)]
 the $8$-element set
$$B(f)=\{(\alpha^3:\alpha:1)\in \P^2(K_a)\mid \alpha\in \RR\}$$
 is in general position.
 \item[(ii)]
 Let $S_{B(f)}$ be the   Del Pezzo surfaces of degree $1$ obtained by blowing up the points of $B(f)$.
 Let $C_ {B(f)}$  be the  branch curve of  $S_{B(f)}$ with respect to its bi-anticanonical map. Let $J(C_ {B(f)})$ be the
  jacobian of $C_ {B(f)}$. Then the endomorphism algebra $\End^0(J(C_ {B(f)}))$ of $J(C_ {B(f)})$
(over $\bar{K}$) is either $\Q$ or a quadratic field; in
particular, $J(C_ {B(f)})$ is an absolutely simple abelian fourfold.
\item[(3)]
Let $h(t)\in K[t]$ be another irreducible  degree $8$ polynomial, whose Galois
group is either $\ST_8$ or $\A_8$.   Assume that the splitting fields of $f$ and $h$ are linearly disjoint over $K$. Then
$J(C_ {B(f)})$ and $J(C_ {B(h)})$     are not isomorphic as abelian varieties over    $\bar{K}$ and therefore the surfaces
$S_{B(f)}$ and $S_{B(h)}$   are not biregularly isomorphic over over $\bar{K}$.
   \end{itemize}
\end{thm}

\begin{ex}
If    $K=\Q$ and
$f_1(x)=x^8-x-1$ then the Galois group of $f_1$ is $\ST_8$ \cite[p. 45, Rem. 2]{SerreG}.
Clearly,  $S_{B(f_1)}$, $C_ {B(f_1)}$   and  $J(C_ {B(f_1)})$ are defined over $\Q$ and, thanks to Theorem \ref{mainex},
$\End^0(J(C_ {B(f_1)}))$ is either $\Q$ or a quadratic field. (Recall  that $J(C_ {B(f_1)})$ has the same endomorphism  algebra over
$\bar{\Q}$ and $\C$.) This implies that $J(C_ {B(f_1)})$ is  simple, viewed as a complex abelian fourfold.
\end{ex}

\begin{ex}
If $K=\Q$,  $\Q(t)$ and  $\bar{\Q}(t)$    are the fields of rational functions in one variable $t$ over $\Q$ and
 $\bar{\Q}$  respectively then the polynomial $f_t(x)=x^8-x-t\in \Q(t)\subset \bar{\Q}(t)$ has Galois group $\ST_8$
over $\bar{\Q}(t)$ \cite[p. 139]{SerreMW}. Now Hilbert's irreducibility theorem \cite[Sect. 10.1]{SerreMW} implies that there exists an infinite set $N$ of rational numbers
 such that for each $n\in N$ the polynomial $f_n(t)=x^8-x-n$ has Galois group   $\ST_8$   over $\Q$ and the splitting fields of $f_n$
 are linearly disjoint for distinct $n$.    Clearly,  $S_{B(f_n)}$, $C_ {B(f_n)}$   and  $J(C_ {B(f_n)})$ are defined over $\Q$ and,
 thanks to Theorem \ref{mainex},
$\End^0(J(C_ {B(f_n)}))$ is either $\Q$ or a quadratic field.   In addition, abelian varieties   $J(C_ {B(f_n)})$'s are not isomorphic
over $\bar{\Q}$ for distinct $n$.
Since all $\C$-isomorphisms among  $J(C_ {B(f_n)})$'s are defined over $\bar{\Q}$, the {\sl complex} abelian varieties $J(C_ {B(f_n)})$'s
are {\sl not} isomorphic for distinct $n$.  This implies that the surfaces    $S_{B(f_n)}$'s are not biregularly isomorphic over $\C$ for distinct $n$.
\end{ex}

The Del Pezzo surfaces involved look rather special. That is why
we prove that the assertion about the endomorphism algebra of the
corresponding jacobian remains true  when the corresponding Galois
image in $W(E_8)$  contains a subgroup that is a conjugate of the
$\A_8$. In particular, the jacobian of the branch curve of
``generic" Del Pezzo surface of degree $1$ is absolutely simple.

On the other hand, assuming that the Galois action on the Picard
group of  a Del Pezzo surface of degree $1$ is {\sl maximal}
(i.e., the Galois image coincides with $W(E_8)$), we prove that
the jacobian of the corresponding branch curve has no non-trivial
endomorphisms. It would be interesting to find explicit examples
of degree $1$ surfaces with maximal Galois action. (See \cite{Ek} for explicit examples
 of cubic surfaces with maximal Galois image.  The case of degree $2$ is discussed in \cite{Erne}).

 The paper is organized as follows. In Section \ref{DP} we discuss
 interrelations between Picard groups of a Del Pezzo surface of
 small degree and the corresponding branch curve. Our
 exposition is based on letters of Igor Dolgachev to the author.
  Section \ref{AV} deals with abelian fourfolds whose Galois
  module of points of order $2$ has a rather special structure.
  Our main results are stated and proved in Section \ref{jacob}.
 In Section \ref{explicit} we describe an algorithm for finding an (explicit) equation for a (singular) plane
  birational model of $C_ {B(f)}$ in terms of $f$.

I am deeply grateful to Igor Dolgachev for his interest to this
paper and generous help.

This work was started during the special semester ``Rational and
integral points on higher-dimensional varieties" at the MSRI (Spring 2006). The
author is grateful to the MSRI and the organizers of this program.
My special thanks go to the referee, whose comments helped to improve the exposition.

\section{Del Pezzo surfaces of degree $1$}
\label{DP}

\begin{sect}
\label{weyl}
 Let $d=1$ or $2$. We write $I^{1,9-d}$ for the  standard
odd unimodular hyperbolic lattice of rank $10-d$ and signature
$(1,9-d)$. This means that $I^{1,9-d}$ is a free $\Z$-module of
rank $10-d$ provided with the unimodular symmetric bilinear form
$( , ):I^{1,9-d} \times I^{1,9-d} \to \Z$ and the standard
orthogonal basis $\{e_0,e_1,\ldots,e_{9-d}\}$ such that
$$(e_0,e_0)=1, \ (e_i,e_i)=-1  \ \forall i\ge 1.$$
We write $\OO(I^{1,9-d})$ for the group of isometries of
$I^{1,9-d}$.

Let us put $\omega_{9-d}:=-3e_0+e_1+\ldots+e_{9-d}$. Clearly,
$(\omega_{9-d},\omega_{9-d})=d \ne 0$. Let us consider the
orthogonal complement $\omega_{9-d}^{\perp}$ of $\omega_{9-d}$ in
$I^{1,9-d}$. Clearly, $\omega_{9-d}^{\perp}$ is a free
$\Z$-(sub)module of rank $9-d$ and the restriction
$$( , ):\omega_{9-d}^{\perp}\times \omega_{9-d}^{\perp} \to \Z$$
is a negative-definite non-degenerate symmetric bilinear form. (If
$d=1$ then it is even unimodular.) It is known \cite[Th.
23.9]{Manin} that there exists an isometry
$$\omega_{9-d}^{\perp} \to  Q(E_{9-d}),$$
where $Q(E_{9-d})$ is the root lattice of type $E_{9-d}$ equipped
with the scalar product with opposite sign. In addition, if we
identify (via this isometry)  $\omega_{9-d}^{\perp}$ and
$Q(E_{9-d})$ then the orthogonal group $\OO(\omega_{9-d}^{\perp})$
of $\omega_{9-d}^{\perp}$ coincides with the corresponding Weyl
group $W(E_{9-d})$ \cite[Th. 23.9]{Manin}. This implies that the
special orthogonal group $\SO(\omega_{9-d}^{\perp})$ of
$\omega_{9-d}^{\perp}$ coincides with the (index $2$ sub)group
$W^{+}(E_{9-d})$ of elements of determinant $1$ in $W(E_{9-d})$.

There is a natural injective homomorphism
$\mu:\ST_{9-d}\hookrightarrow \OO(I^{1,9-d})$ defined as follows.
$$\tau\mapsto \mu(\tau):I^{1,9-d} \to I^{1,9-d}, \ e_0\mapsto e_0, e_i\mapsto
e_{\tau(i)} \ \forall i\ge 1 \eqno(0).$$ This provides $I^{1,9-d}$
with the natural structure of faithful $\ST_{9-d}$-module.
 Clearly, each
$\mu(\tau)$ leaves invariant $\omega_{9-d}$ and therefore induces
an isometry of $\omega_{9-d}^{\perp}$, which we denote by
$$\iota(\tau)\in \OO(\omega_{9-d}^{\perp})=W(E_{9-d}).$$
 Obviously, if $v:I^{1,9-d} \to I^{1,9-d}$ is an isometry that
leaves $\omega_{9-d}$ invariant and coincides with $\iota(\tau)$
on $\omega_{9-d}^{\perp}$ then $v=\mu(\tau)$. This implies that
$\omega_{9-d}^{\perp}$ is a faithful $\ST_{9-d}$-submodule. We
write
$$\iota:\ST_{9-d} \hookrightarrow \OO(\omega_{9-d}^{\perp})=W(E_{9-d}), \ \tau\mapsto \iota(\tau)$$
for the corresponding defining homomorphism.

Suppose now that $d=1$. Then $9-d=8$ and $(\omega_8,\omega_8)=1$.
This gives us a $\ST_8$-invariant orthogonal splitting
$$I^{1,8}=\Z\cdot \omega_8\oplus \omega_{8}^{\perp}$$
and allows us to identify $W(E_{8})$ with a certain subgroup of
$\OO(I^{1,8})$, namely,
$$W(E_{8})=\OO(\omega_{8}^{\perp})=\{v\in \OO(I^{1,8})\mid
v(\omega_8)=\omega_8\}.$$ (Notice that under this identification
$\iota(\tau)\in \OO(\omega_{8}^{\perp})$ goes into $\mu(\tau)\in
\OO(I^{1,8})$.)

\end{sect}

\begin{sect} \label{mark}
 Let $S$ be a Del Pezzo surface over $\bar{K}$ of
degree $d = 1$ (or $2$). It is well-known \cite[Ch. IV, Sect.
24]{Manin} that $\Pic(S)$ is a free $\Z$-module of rank $10-d$
provided with unimodular bilinear intersection form
$$( , ): \Pic(S)\times \Pic(S) \to \Z$$ of signature $(1,9-d)$;
The anticanonical class $-\K_S$ is ample, $(\K_S,\K_S)=d$,
$$\dim_{\bar{K}}H^0(S,\Oc_S(-\K_S))=d+1$$ and
$$H^0(S,\Oc_S(n \K_S))=\{0\}, \ H^1(S,\Oc_S(\pm n \K_S))=\{0\} \eqno(1)$$ for all positive
integers  $n$ \cite[Cor. 3 on p. 65]{Demazure}. If $E$ is an
exceptional curve of the first kind on $S$ then
$$<E\cdot K_S>=-1$$
\cite[Sect. 26]{Manin}. Let us put
$$L: = \K_S^\perp\subset \Pic(S),$$
i.e., $L$ is the orthogonal complement of $\K_S$ in $\Pic(S)$ with
respect to the intersection pairing. Clearly, $L$ is a free
$\Z$-module of rank $9-d$ and (the restriction)
$$( , ): L \times L \to \Z$$
is a non-degenerate symmetric bilinear form. We write $\OO(L)$
(resp. $\SO(L)$) for the group of automorphisms of $L$ (resp.
automorphisms of $L$ with determinant $1$) preserving the
intersection pairing.

 Recall \cite{Dolgachev} that a {\sl marking} of $S$ is an
 isometry
$\phi:\Pic(S) \to I^{1,9-d}$ such that $\phi(\K_S) =\omega_{9-d}$.
It is  known  that a marking always exists \cite[Prop.
25.1]{Manin}. In addition, each marking is a {\sl geometric
marking}, i.e. can be realized in such a way that $\phi^{-1}(e_i)$
(for positive $i$) are the classes of the exceptional curves under
some blow-up $f:S\to \P^2$ \cite[Ch. 8, Sect. 8.2]{Dolgachev2}. A
marking induces an isometry of lattices
$$\phi:L = \K_S^\perp \cong \omega_{9-d}^{\perp}= Q(E_{9-d})$$
and gives rise to a group isomorphism
$$\OO(L) \to W(E_{9-d}), u \mapsto \phi u \phi^{-1}.$$

Let $w_0$ be the nontrivial center element of $W(E_{9-d})$ which
acts as -1 on the root lattice. There exists a unique automorphism
$g_0$ of $S$ (called the Geiser ($d = 2$), or Bertini ($d = 1$)
involution) such that $\phi g_0 \phi^{-1} = w_0$
 \cite[Th. 4.7]{Manin2},
\cite[pp. 66-69]{Demazure}, \cite[Ch. 8, Sect. 8.2]{Dolgachev2}.
Let $S^{g_0}$ be the fixed locus of $g_0$. It is a smooth
irreducible projective curve $C$ of  genus $3$ if $d=2$. If $d=1$ then
$S^{g_0}$ is a disjoint union of an isolated point and a a smooth
irreducible projective curve $C$ of  genus 
 $4$.

We write $J(C)$ for the jacobian of $C$: it is an abelian variety
over $\bar{K}$ of dimension $3$ (resp. $4$). We write $\Pic(C)_2$
for the kernel of multiplication by $2$ in $\Pic(C)$. Clearly,
$\Pic(C)_2$ coincides with the group $J(C)_2$ of points of order
$2$ on $J$; it is a $\F_2$-vector space of dimension $6$ (resp.
$8$) provided with the alternating nondegenerate bilinear form
called the {\sl Weil pairing} \cite{MumfordAV,MumfordT}
$$< , >: J(C)_2 \times J(C)_2 \to \F_2.$$

Recall that if  $D$ is a divisor on $S$ then
$$H^0(S,\Oc_S(D))=\L(S,D):=\{u\in \bar{K}(S)\mid \div(u)+D\ge 0\}\subset
\bar{K}(S).$$
\end{sect}

\begin{lem}
\label{lem0} Let $E$ be an exceptional curve of the first kind on
$S$ and $n$ a positive integer. Then:
\begin{itemize}
 \item[(i)]
  $\L(S,E)=\bar{K},\ \L(S,E+n\K_S)=\{0\}$.
 \item[(ii)]
$H^1(S, \Oc_S(-E+n \K_S))=\{0\}, \ H^1(S, \Oc_S(E-(n-1)
\K_S)=\{0\}$.
\end{itemize}
\end{lem}

\begin{proof}
(i) Since the self-intersection index of irreducible curve $E$ is
negative, the linear system $\mid E\mid$ consists of single
divisor $E$. This implies that $\L(S,E)=\bar{K}$. Clearly,
$$\L(S,E+n\K_S)\subset \L(S,E)=\bar{K}.$$
Since $-n\K_S$ is ample, $\L(S,E+n\K_S)=\{0\}$.

In order to prove (ii), let us consider the exact sequence
$$0\to  \Oc_S(-E+n \K_S)\to  \Oc_S(n \K_S)  \to \Oc_E(n \K_S)\to 0,$$
which leads to the cohomological exact sequence
$$H^0(E, n \K_S\mid E)=H^0(S,\Oc_E(n \K_S))\to H^1(S, \Oc_S(-E+n \K_S)) \to
H^1(S,\Oc_S(n\K_S)).$$ Since $(n \K_S, E)=-n<0$, we have $H^0(E, n
\K_S\mid E)=\{0\}$. By (1), $H^1(S,\Oc_S(nK_S))=\{0\}$. Now the
exactness of the cohomological sequence implies that $H^1(S,
\Oc_S(-E+n \K_S)=\{0\}$. By Serre's duality, we obtain that
$H^1(S, \Oc_S(E-(n-1) \K_S))=\{0\}$.
\end{proof}

Recall \cite{MumfordT} that a divisor class $\eta$ on a smooth
curve $C$ is called a {\sl theta characteristic} if $2\eta =
\K_C$. It is called {\sl even} (odd) if $h^0(\eta)$ is even (odd).
A theta characteristic $\eta$ defines a quadratic form on
$\Pic(C)_2$ by
$$q_\eta(x) = h^0(\eta+x)+h^0(\eta) \bmod 2.$$ The associated
bilinear form $e(x,y) = q_\eta(x+y)+q_\eta(x)+q_\eta(y)$ coincides
with $< , >$ on $J(C)_2=\Pic(C)_2$. In particular, if $\eta$ is an
{\sl even} theta characteristic then
$$<x ,y >=h^0(\eta+x+y)+h^0(\eta+x)+h^0(\eta+y) \bmod 2 \eqno(2).$$

\begin{lem}
\label{lem1}
 Let $r:\Pic(S)\to \Pic(C)$ be the restriction
homomorphism. Then

\begin{itemize}
\item[(i)] $r(a)\in \Pic(C)_2$ for all $a\in L$. In other words,
$$r(2L)=0.$$
 \item[(ii] If $E$ is an exceptional curve of the first kind on
 $S$ then
$r(E)$ is an odd theta characteristic on $C$; more precisely,
$$h^0(C,r(E)):=\dim_{\bar{K}}(H^0(C,r(E)))=1.$$
In addition,  if $d=1$ then $r(-\K_S)$ is a theta characteristic
on $C$.
\end{itemize}
\end{lem}

\begin{proof} (i) Since $g_0$ acts on $C$ as identity map, we have $r(g_0(a)) = g_0(r(a)) =
r(a)$. Since $g_0|L = -1$, we get $r(-a) = r(a)$ and hence $2r(a)
= 0$.

(ii) Recall that if $D$ is a divisor on $C$ then
$$H^0(C,\Oc_C(D))=\L(C,D):=\{u\in \bar{K}(C)\mid \div(u)+D\ge 0\}\subset
\bar{K}(C).$$ If $u$ is a non-constant rational function on $C$
then the degree of the divisor of poles of $u$ coincides with the
degree of field extension $\bar{K}(C)/\bar{K}(u)$ \cite[Th.
2.2]{Lang}.

($d = 2$) Recall that $|-\K_S|$ defines a finite map $\f0$ of degree
2 from $S$ to $\P^2$ and the corresponding involution of $S$ is
$g_0$. Thus the image of $C$ is a smooth plane quartic, and hence
$C\in |-2\K_S|$  \cite[p. 67 ]{Demazure}. For any irreducible
curve $R$ on $S$ we have $R+g_0(R) = \f0^{*}(k\ell)$, where $\ell$
is the class of a line on $P^2$ and $k$ is an integer.
Intersecting with $-\K_S = \f0^{*}(\ell)$ we get $-2\K_S\cdot R =
2k.$ In particular, if $R = E$ is an exceptional curve of the
first kind, then $E\cdot \K_S =-1$, we  get $E+g_0(E) = \f0^*(\ell)
= -\K_S$. Since $E|C = g_0(E)|C$, we obtain $2r(E) = -r(\K_S)$. By
the adjunction formula, $\K_C = r(\K_S+C) = -r(\K_S)$. This proves
that  $r(E)$ is a theta characteristic. It is certainly effective,
and therefore we may view $r(E)$ as the linear equivalence class
of an effective divisor $D$ of degree $2$ on $C$. We need to prove
that $\dim_{\bar{K}}(\L(C,D))=1$.

If $u \in \L(C,D)$
 is  not a constant then the degree of its polar divisor is
 either $1$ or $2$. In the former case, $\bar{K}(C)=\bar{K}(u)$ and $C$
is rational, which is not the case, because $C$ has genus $3$. In
the latter case, $\bar{K}(C)/\bar{K}(u)$ is a quadratic extension
and therefore $C$ is hyperelliptic, which is not the case, because
$C$ is a plane smooth quartic \cite[Ch. 5, Exercise 3.2 and
Example 5.2.1]{H}.
  This implies that $\L(C,D)$ consists of scalars
  and therefore has dimension $1$.

 ($d = 1$) The linear system $|-2\K_S|$ defines a double cover $\f0:S\to Q$,
 where $Q$ is a quadratic  cone in $\P^3$. The branch curve $W$ is a curve of genus 4,
 the intersection of $Q$ with a cubic surface \cite[Ch. 8, Example 8.2.4]{Dolgachev2}.
 This implies that $\f0^*(W) = 2C\in |-6\K_S|$,
 hence $C\in |-3\K_S|$ \cite[pp. 68--69 ]{Demazure}. By adjunction formula,
$\K_C = r(-3\K_S+\K_S) = -2r(\K_S)$. Now similar to the case $d =
2$, we obtain $E+g_0(E) = -2K_S$ and $2r(E) = -2r(\K_S) = \K_C$.
This proves that $r(-\K_S)$ and $r(E)$ are  theta characteristics.
  Clearly, $r(E)$ is  effective and therefore $H^0(C,r(E)) \ne \{0\}$.
The short exact sequence
$$0 \to \Oc_S(E+3\K_S) = \Oc_S(E-C) \to \Oc_S(E) \to \Oc_C(E) \to 0$$
gives rise to the  exact cohomological sequence
$$H^0(S, \Oc_S(E))\to H^0(S,\Oc_C(E))=H^0(C,r(E)) \to
H^1(S,\Oc_S(E+3\K_S)).$$

Suppose that we know that $H^1(S,\Oc_S(E+3\K_S))=\{0\}$. This
implies that $H^0(S, \Oc_S(E))\to H^0(C,r(E))$ is a surjection.
Since $H^0(S, \Oc_S(E))$ is one-dimensional, thanks to Lemma
\ref{lem0}(i) and $H^0(C,r(E)) \ne \{0\}$, we conclude that
$H^0(S, \Oc_S(E))\cong H^0(C,r(E))$; in particular, $h^0(C,r(E)) =
1$.

So, in order to finish the proof, we need to check that
$$H^1(S,\Oc_S(E+3\K_S))=\{0\}.$$
By Serre's duality, it suffices to prove that
$H^1(S,\Oc_S(-E-2\K_S))=\{0\}$. In order to do that, recall that
$-2K_S=E+g_0(E)$ in $\Pic(S)$ where $g_0(E)$ is also an
exceptional curve of the first kind. This implies that
$$-E-2\K_S=-E+E+g_0(E)=g_0(E).$$
Applying Lemma \ref{lem0} to $g_0(E)$ and $n=1$, we conclude that
$H^1(S,\Oc_S(g_0(E)))=\{0\}$ and therefore
$H^1(S,\Oc_S(-E-2\K_S))=H^1(S, \Oc_S(g_0(E)))=\{0\}$.
\end{proof}

\begin{rem}
If $d=1$ then the only one theta characteristic of $C$ with $h^0 >
1$ is the vanishing one equal to $r(-\K_S)$ (see below).
\end{rem}

\begin{rem}
If $d=1$ then $K_C$ is very ample \cite[Ch. 8, Example
8.2.4]{Dolgachev2}. This implies that $C$ is {\sl not}
hyperelliptic.
\end{rem}

\begin{lem}
\label{lem2} Assume $d = 1$. Then:
\begin{itemize}
\item $r(-\K_S)$ is an even theta characteristic on $C$. \item
 Let $\phi:\Pic(S)\to I^{1,8}$ be a marking and
$E_i = \phi^{-1}(e_i), i = 1,\ldots,8$. Then $v_i = E_i+\K_S\in L$
and points $x_i = r(v_i)\in \Pic(C)_2=J(C)_2$ satisfy $<x_i,x_j> =
1$ if $i\ne j$.
\end{itemize}
\end{lem}

\begin{proof}  Let $\eta_0 = r(-\K_S)$. By Lemma \ref{lem1}(ii),
$\eta_0$ is a theta
characteristic. The short exact sequence
$$0 \to \Oc_S(2\K_S) = \Oc_S(-\K_S-C) \to \Oc_S(-\K_S) \to \Oc_C(-\K_S) \to 0$$
gives rise to the  exact cohomological sequence
$$H^0(S,\Oc_S(2\K_S))\to H^0(S,\Oc_S(-\K_S)) \to H^0(C,-\K_S\mid C)=H^0(C,\eta_0) \to
H^1(S,\Oc_S(2\K_S)).$$ Since  $H^0(S,\Oc_S(2\K_S))=0$ and
$H^1(S,\Oc_S(2\K_S))=0$, the $\bar{K}$-vector spaces
$H^0(C,\eta_0)$ and $H^0(S,\Oc_S(-\K_S))$ are isomorphic; in
particular,
  $h^0(C, \eta_0) =
h^0(-\K_S) = 2.$ Thus $\eta_0$ is an even theta characteristic.

We have
$$\eta_0+x_i=r(-\K_S)+r(E_i+\K_S)=r(E_i), \ \eta_0+x_j=r(E_j),$$
$$\eta_0+x_i+x_j=r(E_i)+r(E_j+\K_S)=r(E_i+E_j+\K_S).$$
Applying  (2) to $\eta=\eta_0$, we conclude that
$$<x_i,x_j>= h^0(\eta_0+x_i+x_j)+h^0(\eta_0+x_i)+h^0(\eta_0+x_j) \bmod 2 = $$
$$= h^0(r(E_i+E_j+\K_S))+h^0(r(E_i))+h^0(r(E_j)) \bmod 2=$$
$$=h^0(r(E_i+E_j+\K_S))+1+1 \bmod 2 = h^0(r(E_i+E_j+\K_S))\bmod 2.$$ (Here we
used Lemma \ref{lem1} that tells us that $h^0(r(E_i)) =
h^0(r(E_j)) = 1$.)

Since $r(2L) = 0$ and  $E_j+\K_S \in L$, we obtain that
$r(E_i+E_j+\K_S)=r(E_i-E_j-\K_S)$ and therefore
$$<x_i,x_j> = h^0(r(E_i-E_j-\K_S)) \bmod 2.$$
 Now $$(E_i-E_j-\K_S, \K_S) = -1, \ (E_i-E_j-\K_S, E_i-E_j-\K_S) = -1.$$
 It is known \cite[Th. 26.2(i)]{Manin} that this implies that
$E_i-E_j-\K_S$ is linearly equivalent to the class of an
exceptional curve of the first kind. By Lemma \ref{lem1}, we
obtain $h^0(r(E_i-E_j-\K_S)) = 1$. So, $<x_i,x_j> = 1$.
\end{proof}

Now we need the following elementary result from linear algebra.

\begin{lem}
\label{linalg} Let $\F$ be a field of characteristic $2$, $V$ a
finite-dimensional $\F$-vector space,
$$\phi: V\times V \to \F$$
an alternating $\F$-bilinear form, $c$ a non-zero element of $\F$.
Let $m$ be a positive even integer and $\{z_1, \dots z_m\}$ an
$m$-tuple of vectors in $V$ such that
$$\phi(z_i,z_j)=c \
\forall \ i\ne j.$$

Then:
\begin{itemize}
\item[(i)]
 $\{z_1, \dots z_m\}$ is a set of linearly independent
vectors; in particular, it is a basis if $m=\dim(V)$.

\item[(ii)] Let $V_m$ be the subspace of $V$ generated by all
$z_i$'s. Then the restriction of $\phi$ to $V_m$ is nondegenerate.
In particular, if $m=\dim(V)$ then $V_m=V$ and $\phi$ is
non-degenerate.
\end{itemize}
\end{lem}

Both assertions of Lemma \ref{linalg} follow immediately from the next statement.

\begin{prop}
\label{vanish}
 Let $\{a_1\dots , a_m\}$ be an $m$-tuple of elements of $\F$
such that $z=\sum_{i=1}^m a_i z_i$ satisfies $\phi(z, z_j)=0 \
\forall \  j$. Then all $a_i$'s vanish.
\end{prop}

\begin{proof}[Proof of Proposition \ref{vanish}] Let us put
$$b=\sum_{i=1}^m a_i.$$
Since for all $j$ with $1\le j\le r$
$$\phi(z_j,z_j)=0= \phi(z,z_j)=\phi\left(\sum_{i=1}^r a_i z_i,z_j\right)=\sum_{i\ne j} a_i\cdot c,$$
we have $c\cdot\sum_{i\ne j} a_i=0$. This implies that
 $b-a_j=\sum_{i\ne j} a_i=0 \ \forall j$.
 It follows that
 $$0=\sum_{j=1}^m(b-a_j)=mb-\sum_{j=1}^m a_j=mb-b=(m-1)b=b,$$
 since $m$ is even and $\fchar(\F)=2$. So, $b=0$. Since every $b-a_j=0$, we conclude that every $a_j=0$.
\end{proof}

\begin{thm}
\label{CS2}
 Assume $d = 1$. The map $L/2L\to \Pic(C)_2$ induced by $r$ is an
isometry with respect to the bilinear intersection  form $< , >$
on $L$ reduced modulo 2 and the Weil pairing $< , >$ on
$\Pic(C)_2$.
\end{thm}

\begin{proof} Let $v_1,\ldots,v_8\in L$ be as in Lemma \ref{lem2}. We have
$$(v_i,v_j)=(E_i+\K_S,E_j+\K_S) = -1$$ and therefore $1=(v_i,v_j)\bmod
2= <r(v_i),r(v_j)>$.  It follows from Lemma \ref{linalg} that
$r(v_i)$'s form a basis in $\Pic(C)_2$.
This implies that
$\{r(v_1),\ldots,r(v_8)\}$ is a basis in $L/2L$ and moreover the
restriction map $r \bmod 2$ of the 8-dimensional vector spaces
over $\F_2$ preserves the corresponding non-degenerate bilinear
forms. Obviously it implies that  the map is an isometry.
\end{proof}

\begin{rem} A similar assertion is true in the case $d = 2$ and
is proven in \cite[pp. 160--162]{Dolgachev}. One can give it a
similar proof without using the Smith exact sequence. We put $x_i
= r(E_i-E_1)$, where $\{E_1,\ldots,E_7\}$ are defined similar to
the above ($i=2,\dots ,7$). We take the odd theta characteristic
$\eta = r(E_1)$ and compute the associated bilinear form of the
quadratic form $q_\eta$. We do the similar computation, using the
fact that the seven odd theta characteristics $\eta_i = r(E_i)$
form an {\sl Aronhold set}, i.e., $\eta_i+\eta_j-\eta_k$ is an
even theta characteristic for any distinct $i,j,k$ \cite[pp.
166--169]{Dolgachev}. (See also \cite[Ch. 6, Sect.
6.1]{Dolgachev2}.)
\end{rem}

\begin{rem}
\label{GCS2} Suppose that $S$ is defined over a field $K$. Then
 $\Pic(S)$ carries the natural structure of $\Gal(K)$-module, the intersection pairing and
$\K_S$ are Galois-invariant \cite[Ch. IV]{Manin}. In addition,
$g_0$ is defined over $K$ \cite[Th. 4.7]{Manin2} and therefore $C$
is a $K$-curve on $S$. It follows easily from Theorem \ref{CS2}
that $L$ is a Galois submodule and $L/2L\to \Pic(C)_2$ is an
isomorphism of Galois modules.
\end{rem}

\section{Abelian varieties}
\label{AV}
 A surjective homomorphism of finite groups
$\pi:\GG_1\twoheadrightarrow \GG$ is called a {\sl minimal cover}
if no proper subgroup of $\GG_1$ maps onto $\GG$ \cite{FT}.
Clearly, if $\GG$ is perfect and $\GG_1\twoheadrightarrow \GG$ is a
minimal cover then $\GG_1$ is also perfect.

Let $F$ be a field, $F_a$ its algebraic closure  and
$\Gal(F):=\Aut(F_a/F)$ the absolute Galois group of $F$. If $X$ is
an abelian variety of positive dimension
 over $F_a$ then we write $\End(X)$ for the ring of all its
$F_a$-endomorphisms and $\End^0(X)$ for the corresponding
$\Q$-algebra $\End(X)\otimes\Q$. We write $\End_F(X)$ for the ring
of all  $F$-endomorphisms of $X$ and $\End_F^0(X)$ for the
corresponding  $\Q$-algebra $\End_F(X)\otimes\Q$ and $\CC$ for the center of $\End^0(X)$. Both $\End^0(X)$
and $\End_F^0(X)$  are semisimple finite-dimensional
$\Q$-algebras.

 The  group
$\Gal(F)$ of $F$ acts on $\End(X)$ (and therefore on $\End^0(X)$)
by ring (resp. algebra) automorphisms and
$$\End_F(X)=\End(X)^{\Gal(F)}, \ \End_F^0(X)=\End^0(X)^{\Gal(F)},$$
since every endomorphism of $X$ is defined over a finite separable
extension of $F$.

If $n$ is a positive integer that is not divisible by $\fchar(F)$
then we write $X_n$ for the kernel of multiplication by $n$ in
$X(F_a)$; the commutative group  $X_n$ is a free $\Z/n\Z$-module
of rank $2\dim(X)$ \cite{MumfordAV}. In particular, if $n=2$ then
$X_{2}$ is an $\F_{2}$-vector space of dimension $2\dim(X)$.

 If $X$ is defined over $F$ then $X_n$ is a Galois
submodule in $X(F_a)$ and all points of $X_n$ are defined over a
finite separable extension of $F$. We write
$\bar{\rho}_{n,X,F}:\Gal(F)\to \Aut_{\Z/n\Z}(X_n)$ for the
corresponding homomorphism defining the structure of the Galois
module on $X_n$,
$$\tilde{G}_{n,X,F}\subset
\Aut_{\Z/n\Z}(X_{n})$$ for its image $\bar{\rho}_{n,X,F}(\Gal(F))$
and $F(X_n)$ for the field of definition of all points of $X_n$.
Clearly, $F(X_n)$ is a finite Galois extension of $F$ with Galois
group $\Gal(F(X_n)/F)=\tilde{G}_{n,X,F}$. If $n=2$  then we get a
natural faithful linear representation
$$\tilde{G}_{2,X,F}\subset \Aut_{\F_{2}}(X_{2})$$
of $\tilde{G}_{2,X,F}$ in the $\F_{2}$-vector space $X_{2}$.

Now and till the end of this Section we assume  that $\fchar(F)\ne
2$. It is known \cite{Silverberg} that all endomorphisms of $X$
are defined over $F(X_4)$; this gives rise to the natural
homomorphism
$$\kappa_{X,4}:\tilde{G}_{4,X,F} \to \Aut(\End^0(X))$$ and
$\End_F^0(X)$ coincides with the subalgebra
$\End^0(X)^{\tilde{G}_{4,X,F}}$ of $\tilde{G}_{4,X,F}$-invariants
\cite[Sect. 1]{ZarhinLuminy}.

The field inclusion $F(X_2)\subset F(X_4)$ induces a natural
surjection \cite[Sect. 1]{ZarhinLuminy}
$$\tau_{2,X}:\tilde{G}_{4,X,F}\twoheadrightarrow\tilde{G}_{2,X,F}.$$

The following definition has already appeared  in \cite{ElkinZ}.

\begin{defn} We say that $F$ is 2-{\sl balanced} with respect to
$X$ if $\tau_{2,X}$ is a minimal cover.
\end{defn}

\begin{rem}
\label{overL}
 Clearly, there always exists a subgroup $H
\subset \tilde{G}_{4,X,F}$ such that $H\to\tilde{G}_{2,X,F}$ is
surjective and a minimal cover. Let us put $L=F(X_4)^H$. Clearly,
$$F \subset L \subset F(X_4), \ L\bigcap F(X_2)=F$$
and $L$ is a maximal overfield of $F$ that enjoys these
properties. It is also clear that There exists an   overfield $L$
such that
$$F \subset L \subset F(X_4), \ L\bigcap F(X_2)=F,$$
$$F(X_2)\subset L(X_2),\ L(X_4)=F(X_4),
\tilde{G}_{2,X,L}=\tilde{G}_{2,X,F}$$
 and $L$ is $2$-{\sl balanced}   with respect to
$X$ (see \cite[Remark 2.3]{ElkinZ}).
\end{rem}

\begin{thm}
\label{rep}
 Suppose that $E:=\End_F^0(X)$ is a field that contains the
center $\CC$ of $\End^0(X)$. Let $\CC_{X,F}$ be the centralizer of
$\End_F^0(X)$ in $\End^0(X)$.

Then:

\begin{itemize}
\item[(i)] $\CC_{X,F}$ is a central simple $E$-subalgebra in
$\End^0(X)$. In addition, the centralizer of $\CC_{X,F}$ in
$\End^0(X)$ coincides with $E=\End_F^0(X)$ and
$$\dim_E(\CC_{X,F})=\frac{\dim_{\CC}(\End^0(X))}{[E:\CC]^2}.$$
 \item[(ii)] Assume that $F$ is $2$-balanced with
respect to $X$ and $\tilde{G}_{2,X,F}$ is a non-abelian simple
group. If $\End^0(X)\ne E$ (i.e., not all endomorphisms of $X$ are
defined over $F$) then there exist a finite perfect group $\Pi
\subset \CC_{X,F}^{*}$ and a surjective homomorphism $\Pi \to
\tilde{G}_{2,X,F}$ that is a minimal cover. In addition, the
induced homomorphism $E[\Pi] \to \CC_{X,F}$ is surjective, i.e.,
$\CC_{X,F}$ is isomorphic to a direct summand of the group algebra
$E[\Pi]$.
\end{itemize}
\end{thm}

\begin{proof}
This is Theorem 2.4 of \cite{ElkinZ} and Theorem 3.1 of \cite{ElkinZ2}.
\end{proof}

\begin{lem}
\label{Ksimple} Assume that $X_2$ does not contain  proper
$\tilde{G}_{2,X,F}$-invariant even-dimensional subspaces and the
centralizer $\End_{\tilde{G}_{2,X,F}}(X_2)$ has $\F_2$-dimension
$2$.

Then $X$ is $F$-simple and $\End_F^0(X)$ is either $\Q$ or a
quadratic field.
\end{lem}

\begin{proof}
If $Y$ is a proper abelian $F$-subvariety of $X$ then $Y_2$ is a
proper Galois-invariant $2\dim(Y)$-dimensional subspace in $X_2$.
So, $Y_2$ is even-dimensional and $\tilde{G}_{2,X,F}$-invariant.
This implies that such $Y$ does not exist, i.e., $X$ is
$F$-simple. This implies that $\End_F(X)$ has no zero-divisors.
This implies that $\End_F^0(X)$ is a finite-dimensional
$\Q$-algebra without zero divisors and therefore is division
algebra over $\Q$.
 On the other hand, the action of $\End_F(X)$  on
$X_2$ gives rise to an embedding
$$\End_F(X)\otimes \Z/2\Z\hookrightarrow
\End_{\Gal(F)}(X_2)=\End_{\tilde{G}_{2,X,F}}(X_2).$$ This implies
that the rank of free $\Z$-module $\End_F(X)$ does not exceed $2$,
i.e., is either $1$ or $2$. It follows that $\End_F^0(X)$ has
$\Q$-dimension $1$ or $2$ and therefore is commutative. Since
$\End_F^0(X)$ is division algebra, it is a field. If
$\dim_{\Q}\End_F^0(X)=1$ then $\End_F^0(X)=\Q$. If
$\dim_{\Q}\End_F^0(X)=2$ then $\End_F^0(X)$ is a quadratic field.
\end{proof}

\begin{lem}
\label{centerDim} Let us  assume that $g:=\dim(X)>0$ and the center of
$\End^0(X)$ is a field, i.e, $\End^0(X)$ is a simple $\Q$-algebra.

Then:

\begin{itemize}
 \item[(i)]
$\dim_{\Q}(\End^0(X))$ divides $(2g)^2$. \item[(ii)] If
$\dim_{\Q}(\End^0(X))=(2g)^2$ then $\fchar(F)>0$ and $X$ is a supersingular
abelian variety.
\end{itemize}
\end{lem}

\begin{proof}
(ii) is proven in \cite[Lemma 3.1]{ZarhinMRL} (even without any
assumptions on the center).

In order to prove (ii), notice that there exist an absolutely
simple abelian variety $Y$ over $F_a$ and a positive integer $r$
such that $X$ is isogenous (over $F_a$) to a self-product $Y^r$.
We have
$$\dim(X)=r\dim(Y), \End^0(X)=\M_r(\End^0(Y)),
\dim_{\Q}(\End^0(X))=r^2 \dim_{\Q}(\End^0(Y)).$$

 It follows readily from Albert's
classification (\cite[Sect. 21]{MumfordAV}, \cite{Oort}) that
$\dim_{\Q}(\End^0(Y))$ divides $(2\dim(Y))^2$. The rest is clear.
\end{proof}

Let $B$ be an $8$-element set. We write $\Perm(B)$ for the group
of all permutations of $B$. The choice of ordering on $B$
establishes an isomorphism between $\Perm(B)$ and the symmetric
group $\mathrm{S}_8$. We write $\Alt(B)$ for the only subgroup of
index $2$ in $\Perm(B)$. Clearly, every isomorphism
$\Perm(B)\cong\mathrm{S}_8$ induces an isomorphism between
$\Alt(B)$ and the alternating group $\A_8$. Let us consider the
$8$-dimensional $\F_2$-vector space $\F_2^B$ of all $\F_2$-valued
functions on $B$ provided with the natural structure of faithful
$\Perm(B)$-module. Notice that the standard symmetric bilinear
form
$$\F_2^B \times \F_2^B\to \F_2, \ \phi,\psi\mapsto \sum_{b\in
B}\phi(b)\psi(b)$$ is non-degenerate and $\Perm(B)$-invariant.

Since $\Alt(B)\subset \Perm(B)$, one may view $\F_2^B$ as faithful
$\Alt(B)$-module.

\begin{lem}
\label{AA8}
\begin{itemize}
\item[(i)] The centralizer $\End_{\Alt(B)}(\F_2^B)$ has
$\F_2$-dimension $2$.

\item[(ii)] $\F_2^B$ does not contain a proper $\Alt(B)$-invariant
even-dimensional subspace.
\end{itemize}
\end{lem}

\begin{proof}
Since $\Alt(B)$ is doubly transitive, (i) follows from \cite[Lemma
7.1]{Passman}.

Notice that the subspace of $\Alt(B)$-invariants
$$M_0:=(\F_2^B)^{\Alt(B)}=\F_2\cdot 1_B,$$
where $1_B$ is the constant function $1$. In order to prove (ii),
recall that
$$M_{0}\subset M_{1}\subset \F_2^B$$
where $M_{1}$ is the hyperplane of functions with zero sum of
values. It is known \cite{Mortimer} that $M_1/M_0$ is a simple
$\Alt(B)$-module; clearly, $\dim(M_1/M_0)=6$. It follows that if
$W$ is a a proper even-dimensional $\Alt(B)$-invariant subspace of
$\F_2^B$ then $\dim(W)=2$ or $6$. Clearly, the orthogonal
complement $W'$ of $W$ in $\F_2^B$ with respect to the standard
bilinear form is also $\Alt(B)$-invariant and
$\dim(W)+\dim(W')=8$. It follows that either $\dim(W)=2$ or
$\dim(W')=2$. On the other hand, $\Alt(B)=\mathrm{A}_8$ must act
trivially on any two-dimensional $\F_2$-vector space, since $\A_8$
is simple non-abelian and $\GL_2(\F_2)$ is solvable. Since the
subspace of $\Alt(B)$-invariants is one-dimensional, we conclude
that there are no two-dimensional $\Alt(B)$-invariant subspaces of
$\F_2^B$. The obtained contradiction proves the desired result.
\end{proof}

\begin{thm}
\label{mainAV} Let $X$ be a four-dimensional abelian variety over
$F$. Suppose that there exists a group isomorphism
$\tilde{G}_{2,X,F}\cong \Alt(B)$ such that the $\Alt(B)$-module
$X_2$ is isomorphic to $\F_2^B$.

Then one of the following two conditions holds:

\begin{itemize}
\item[(i)] $\End^0(X)$ is either $\Q$ or a quadratic field. In
particular, $X$ is absolutely simple. \item[(ii)] $\fchar(F)>0$
and $X$ is a supersingular abelian variety.
\end{itemize}
\end{thm}

\begin{rem}
\label{help} Lemmas \ref{Ksimple} and \ref{AA8} and Remark
\ref{overL} imply that in the course of the proof of Theorem
\ref{mainAV}, we may assume that $\End_F^0(X)$ is either $\Q$ or a
quadratic field and $F$ is 2-{\sl balanced} with respect to $X$;
in particular, we may assume that $\tilde{G}_{4,X,F}$ is perfect,
since $\tilde{G}_{2,X,F}=\A_8$ is perfect.
\end{rem}


 \begin{proof}[Proof of Theorem \ref{mainAV}]
Following Remark \ref{help}, we assume that $\tilde{G}_{4,X,F}$ is
perfect, $\tau_{2,X}:\tilde{G}_{4,X,F}\twoheadrightarrow
\tilde{G}_{2,X,F}=\A_8$ is a minimal cover and $\End_F^0(X)$ is
either $\Q$ or a quadratic field. Since $\tilde{G}_{4,X,F}$ is
perfect, it does not contain a subgroup of index $2,3$ or $4$. Recall that
$\CC$ is the center of $\End^0(X)$.

\begin{lem}
\label{centerK}
 Either $\CC=\Q\subset \End_F^0(X)$ or
$\CC=\End_F^0(X)$ is a quadratic field.
\end{lem}

\begin{proof}[Proof of Lemma \ref{centerK}]
Suppose that $\CC$ is not a field. Then it is a direct sum
$$\CC=\oplus_{i=1}^r \CC_i$$
of number fields $\CC_1, \dots , \CC_r$ with $1<r\le \dim(X)=4$.
Clearly, the center $\CC$ is a $\tilde{G}_{4,X,F}$-invariant
subalgebra of $\End^0(X)$; it is also clear that
$\tilde{G}_{4,X,F}$ permutes summands $\CC_i$'s. Since
$\tilde{G}_{4,X,F}$ does not contain proper subgroups of index
$\le 4$, each $\CC_i$ is $\tilde{G}_{4,X,F}$-invariant. This implies
that the $r$-dimensional $\Q$-subalgebra
$$\oplus_{i=1}^r \Q \subset \oplus_{i=1}^r \CC_i$$
consists of $\tilde{G}_{4,X,F}$-invariants and therefore lies in
 $\End_F^0(X)$. It follows that $\End_F^0(X)$ has zero-divisors,
 which is not the case. The obtained contradiction proves that $\CC$
 is a field.

 It is known \cite[Sect. 21]{MumfordAV} that $\CC$ contains a totally real number (sub)field $\CC_0$ with
 $[\CC_0:\Q]\mid \dim(X)$ and such that either $\CC=\CC_0$ or  $\CC$ is a purely imaginary
 quadratic extension of $C_0$. Since $\dim(X)=4$, the degree
 $[\CC_0:\Q]$ is $1,2$ or $4$; in particular, the order of
 $\Aut(\CC_0)$ does not exceed $4$. Clearly, $\CC_0$ is
 $\tilde{G}_{4,X,F}$-invariant; this gives us the natural
 homomorphism $\tilde{G}_{4,X,F}\to \Aut(\CC_0)$, which must be
 trivial and therefore $\CC_0$ consists of
 $\tilde{G}_{4,X,F}$-invariants. This implies that
 $\tilde{G}_{4,X,F}$ acts on $\CC$ through a certain homomorphism
 $\tilde{G}_{4,X,F}\to\Aut(\CC/\CC_0)$  and this homomorphism is
 trivial, because the order of $\Aut(\CC/\CC_0)$ is either $1$ (if
 $\CC=\CC_0$) or $2$ (if $\CC\ne \CC_0$). So, the whole $\CC$ consists  of
 $\tilde{G}_{4,X,F}$-invariants, i.e.,
  $$\CC\subset \End^0(X)^{\tilde{G}_{4,X,F}}=\End_F^0(X).$$
  This implies that if $\CC\ne\Q$ then $\End_F^0(X)$ is also not $\Q$
  and therefore is a quadratic field containing $\CC$, which
  implies that $\CC=\End_F^0(X)$ is also a quadratic field.
\end{proof}

It follows that $\End^0(X)$ is a simple $\Q$-algebra (and a
central simple $\CC$-algebra). Let us put $E:=\End_F^0(X)$ and
denote by $\CC_{X,F}$ the centralizer of $E$ in $\End^0(X)$. We have
$$\CC\subset E\subset \CC_{X,F}\subset \End^0(X).$$
Combining Lemma \ref{centerK} with Theorem \ref{rep} and Lemma
\ref{centerDim}, we obtain the following assertion.

\begin{prop}
\label{rep1}
\begin{itemize}
\item[(i)] $\CC_{X,F}$ is a central simple $E$-subalgebra in
$\End^0(X)$,
$$\dim_E(\CC_{X,F})=\frac{\dim_{\CC}(\End^0(X))}{[E:\CC]^2}$$
and $\dim_E(\CC_{X,F})$ divides $(2\dim(X))^2=8^2$.
 \item[(ii)] If $\End^0(X)\ne E$ (i.e., not all endomorphisms of $X$ are
defined over $F$) then there exist a finite perfect group $\Pi
\subset \CC_{X,F}^{*}$ and a surjective homomorphism $\pi:\Pi \to
\tilde{G}_{2,X,F}$ that is a minimal cover.
\end{itemize}
\end{prop}

{\bf End of  Proof of Theorem \ref{mainAV}}.
If $\End^0(X)=E$ then we are done. If $\dim_E(\CC_{X,F})=8^2$ then
$\dim_{\Q}(\End^0(X))\ge \dim_{\CC}(\End^0(X))\ge
\dim_E(\CC_{X,F})=8^2=(2\dim(X))^2$ and it follows from Lemma
\ref{centerK} that $\dim_{\Q}(\End^0(X))=(2\dim(X))^2$ and $X$ is
a supersingular abelian variety. So, further we may and will
assume that
$$\End^0(X)\ne E, \ \dim_E(\CC_{X,F}) \ne 8^2.$$
We need to arrive to a contradiction.
 Let $\Pi \subset
\CC_{X,F}^{*}$ be as in \ref{rep1}(ii). Since $\Pi$ is perfect,
$\dim_E(\CC_{X,F})>1$. It follows from Proposition
 \ref{rep1}(i) that $\dim_E(\CC_{X,F})=d^2$ where $d=2$ or $4$.

 Let us fix an embedding $E\hookrightarrow \C$ and an isomorphism
 $\CC_{X,F}\otimes_E\C\cong \M_d(\C)$. This gives us an embedding
 $\Pi\hookrightarrow \GL(d,\C)$. Further we will identify $\Pi$
 with its image in $\GL(d,\C)$. Clearly, only central elements of
 $\Pi$ are scalars. It follows that there is a central subgroup
 $\ZZ$ of $\Pi$ such that the natural homomorphism $\Pi/\ZZ\to
 \PGL(d,\C)$ is an embedding. The simplicity of
$\tilde{G}_{2,X,F}=\A_8$ implies that $Z$ lies in the kernel of
$\Pi \twoheadrightarrow \tilde{G}_{2,X,F}=\A_8$ and the induced
map $\Pi/Z \to \tilde{G}_{2,X,F}$ is also a minimal cover.
However, the smallest possible degree of nontrivial projective
representation of $\A_8\cong {\mathbf L}_4(2)$ (in characteristic
zero) is $7>d$ \cite{Atlas}.
 Applying Theorem on p. 1092 of \cite{FT} and  Theorem 3 on p. 316
of \cite{KL}, we obtained a desired contradiction.
\end{proof}

  \begin{thm}
  \label{nonisom}
  Suppose that $X$ is as in Theorem \ref{mainAV}.   Suppose that $Y$ is an abelian fourfold over $F$ that enjoys
  one of the following properties:

   \begin{itemize}
   \item[(i)]
   $\tilde{G}_{2,Y,F}$ is solvable;
   \item[(ii)]
   The fields $F(X_2)$ and $F(Y_2)$ are linearly disjoint over $F$.
   \end{itemize}

   If $\fchar(F)=0$ then $X$ and $Y$ are not isomorphic over $\bar{K}$.
  \end{thm}

  \begin{proof}[Proof of Theorem \ref{nonisom}]
  Replacing $F$ by $F(Y_2)$, we may and will assume that  $\tilde{G}_{2,Y,F}=\{1\}$, i.e.,
   the Galois module $Y_2$ is trivial. Clearly, the Galois modules $X_2$ and $Y_2$ are {\sl not} isomorphic.
   By Theorem  \ref{mainAV}, $\End^0(X)$ is either $\Q$ or a quadratic field say, $L$. In the former case all
   the endomorphisms of $X$ are defined over $F$. In the latter case, all
   the endomorphisms of $X$ are defined either over $F$ or over a certain quadratic extension of $F$, because the automorphism
    group of $L$ is the cyclic group of order $2$. Replacing if necessary $F$ by the corresponding quadratic extension, we may and will assume
    that all the endomorphisms of $X$ are defined over $F$. In particular, all the automorphisms of $X$ are defined over $F$;
    in both cases every finite subgroup of $\Aut(X)$ is a finite cyclic group,
    because $\Aut(X)\subset \End^0(X)^{*}$ and $\End^0(X)$ is a field.

     Let $u:X \to Y$ be an $\bar{F}$-isomorphism of abelian varieties.  We need to arrive to a contradiction.
    Since  the Galois modules $X_2$ and $Y_2$ are {\sl not} isomorphic,   $u$ is {\sl not} defined over $F$. Let us consider
    the cocycle
    $$c:\Gal(K)\to \Aut(X), \ \sigma \mapsto c_{\sigma}:= u^{-1} \sigma{u}.$$
    Since the Galois group acts trivially on $\Aut(X)$, the map $c: \Gal(K) \to \Aut(X)$ is a (continuous) group homomorphism.
     Since $u$ is defined over a finite Galois extension of $F$, the image of
    $c$ is a finite subgroup of $\Aut(X)$  and therefore is a finite cyclic group.
    This implies that there is a finite cyclic extension $F'/F$ such that
    $$c_{\sigma}=1 \ \forall \sigma \in \Gal(\bar{F}/F')=\Gal(F')\subset \Gal(F).$$
    It follows that $u$ is defined over $F'$ and therefore the  $\Gal(F')$-modules $X_2$ and $Y_2$ are isomorphic. Clearly,
    the $\Gal(F')$-module $Y_2$ remains trivial.
    However, since $F'/F$ is cyclic and  $\tilde{G}_{2,X,F}\cong \A_8$ is simple non-abelian,
    $$\tilde{G}_{2,X,F'}=\tilde{G}_{2,X,F}\cong \A_8$$
    and therefore   the $\Gal(F')$-module $X_2$ is {\sl not} trivial.  This implies that $X_2$ is not isomorphic to $Y_2$ as  $\Gal(F')$-module
     and we get a desired contradiction.
     \end{proof}

\begin{sect}
Let us consider the $8$-element set $\B=\{1,2, \dots , 7,8\}$.
Then $$\Perm(\B)=\ST_8, \Alt(\B)=\A_8, \F_2^{\B}=\F_2^8.$$ Let us
put $d=1$. Then $9-d=8$ and $\omega_{8}^{\perp}$ carries the
natural structure of $\ST_8$-submodule in $I^{1,8}$ (Sect.
\ref{weyl}).   The inclusion  $\A_8\subset  \ST_8$ provides  $\omega_{8}^{\perp}$
with  the natural structure of $\A_8$-module.
\end{sect}

\begin{lem}
\label{F8S} The $\ST_8$-modules $\F_2^8$ and
$\omega_{8}^{\perp}\otimes\Z/2\Z$ are isomorphic.  In particular,
the   $\A_8$-modules $\F_2^8$ and
$\omega_{8}^{\perp}\otimes\Z/2\Z$ are also isomorphic.
\end{lem}

\begin{proof}
Recall that there is an orthogonal $\ST_8$-invariant splitting
$$I^{1,8}=\Z\cdot \omega_{8}\oplus \omega_{8}^{\perp}.$$
Let us consider the nine-dimensional $\F_2$-vector space
$I^{1,8}\otimes\Z/2\Z$. We write $\bar{e}_i$ for the image of
$e_i$ in $I^{1,8}\otimes\Z/2\Z$. The image of $\omega_{8}$ in
$I^{1,8}\otimes\Z/2\Z$ coincides with
$$\bar{\omega}=\sum_{i=0}^8 \bar{e}_i.$$
The unimodular pairing on $I^{1,8}$ induces a non-degenerate
pairing
$$(I^{1,8}\otimes\Z/2\Z)\times (I^{1,8}\otimes\Z/2\Z)\to \F_2;$$
with respect to this pairing,  all $\bar{e}_i$'s constitute an
orthonormal basis of $I^{1,8}\otimes\Z/2\Z$.

 Clearly, the splitting
$$I^{1,8}\otimes\Z/2\Z=\F_2\cdot\bar{\omega}\oplus
(\omega_{8}^{\perp}\otimes\Z/2\Z)$$ is $\ST_8$-invariant and
orthogonal. In particular, $\omega_{8}^{\perp}\otimes\Z/2\Z$
coincides with the orthogonal complement of $\bar{\omega}$ in
$I^{1,8}\otimes\Z/2\Z$.

It follows that
$$\omega_{8}^{\perp}\otimes\Z/2\Z=\{\sum_{i=0}^8 c_i \bar{e}_i\mid
\sum_{i=0}^8 c_i=0\}.$$ Now the map $$\sum_{i=0}^8 c_i \bar{e}_i
\mapsto \{c_i\}_{i=1}^8\in \F_2^8$$ establishes an isomorphism
between the $\ST_8$-modules $\omega_{8}^{\perp}\otimes\Z/2\Z$ and
$\F_2^8$.
\end{proof}


\begin{sect}
\label{Oref}
We refer to \cite[Ch. VI, Sect. 4, Ex. 1]{Bourbaki}, \cite[p.
85]{Atlas} for the definition and basic properties of the finite
simple non-abelian group $\OO^{+}(8)=\OO^{+}(q_8)$ of order
$2^{12}\cdot 3^5\cdot 5^2\cdot 7$.   In particular, every maximal subgroup in
$\OO^{+}(8)$ has index $\ge 120$ \cite[p. 85]{Atlas}.  This implies that every subgroup of
           $\OO^{+}(8)$ (except $\OO^{+}(8)$ itself) has index $>8$.

  It is known \cite[p. 232]{AtlasB} that every faithful
representation of $\OO^{+}(8)$ in characteristic $2$ has dimension
$\ge 8$. It follows that every faithful  representation of
$\OO^{+}(8)$ in a $8$-dimension $\F_2$-vector space  is
absolutely irreducible and does not split into a tensor product of
two non-trivial representations.
\end{sect}

We refer to \cite{ZarhinTexel,ZarhinVery} for the definition and basic properties of {\sl very simple}
representations.

\begin{thm}
\label{AVO8} Let $X$ be a four-dimensional abelian variety over
$F$. Suppose that there exists a group isomorphism
$\tilde{G}_{2,X,F}\cong \OO^{+}(8)$.

Then one of the following two conditions holds:

\begin{itemize}
\item[(i)] $\End(X)=\Z$.  In particular, $X$ is absolutely simple.
\item[(ii)] $\fchar(F)>0$ and $X$ is a supersingular abelian
variety.
\end{itemize}
\end{thm}

\begin{proof}
Since the natural representation of
$\OO^{+}(8)$ in the $8$-dimension $\F_2$-vector space  $X_2$ is faithful,   it is
absolutely irreducible and does not split into a tensor product of
two non-trivial representations, thanks to Sect. \ref{Oref}. Since every  subgroup in
$\OO^{+}(8)$ (except $\OO^{+}(8)$ itself) has index $>8$ (\ref{Oref}), this
($8$-dimensional) representation is not induced from a subgroup.
It follows from \cite[Th. 7.7]{ZarhinTexel} that the
$\OO^{+}(8)$-module $X_2$ is {\sl very simple}. Now our
 Theorem follows from \cite[Lemma 2.3]{ZarhinTexel}.
\end{proof}

\section{Jacobians}
\label{jacob}

\begin{sect}
\label{general}

Let $K$ be a field of characteristic zero, $\bar{K}$ its algebraic
closure and $\Gal(K):=\Aut(\bar{K}/K)$ the absolute Galois group
of $K$.

 We  use the notation of Section \ref{DP}.
Let $S$ be a Del Pezzo surface of degree $d=1$  over $\bar{K}$,
let $L$ be the orthogonal complement of $\K_S$ in $\Pic(S)$ with
respect to the intersection pairing. Let us fix a marking
$\phi:\Pic(S) \to I^{1,8}$.

Suppose that $S$ is defined over $K$. Recall (Remark \ref{GCS2})
that the Bertini involution on $S$ and the corresponding branch
curve $C$ are also defined over $K$. In addition, there is the
natural Galois action of $\Gal(K)$ on $\Pic(S)$, which preserves
the intersection pairing and $\K_S$. This action is defined by a
certain continuous homomorphism
$$\rho_S:\Gal(K) \to \Aut(\Pic(S)),$$
whose (finite) image consists of isometries; in addition,
$$\rho_S(\sigma)(\K_S)=\K_S \ \forall \sigma\in \Gal(K).$$
The continuous homomorphism
$$\rho_{\phi}:\Gal(K) \to \OO(I^{1,8}), \ \sigma\mapsto
\phi\rho_S(\sigma)\phi^{-1}$$ provides $I^{1,8}$ with a structure
of $\Gal(K)$-module. Clearly, $\phi:\Pic(S) \to I^{1,8}$ is an
isomorphism (and isometry) of Galois modules.

Since $\phi(\K_S)=\omega_8$, the Galois group leaves invariant
$\omega_8$. This implies that $\phi(L)=\omega_8^{\perp}$ and
$$\phi:L \to \omega_8^{\perp} $$
is an isomorphism (isometry) of Galois modules and therefore the
Galois modules $L/2L$ and $\omega_8^{\perp}\otimes \Z/2Z$ are
isomorphic. In addition,
$$I^{1,8}=\Z\cdot \omega_8\oplus \omega_8^{\perp}$$
is a Galois-invariant orthogonal splitting and therefore
$$\rho_{\phi}(\Gal(K))\subset \OO(\omega_8^{\perp})\subset
\OO(I^{1,8}).$$ Here we identify $\OO(\omega_8^{\perp})$ with
 the subgroup of $\OO(I^{1,8})$ that consists of all isometries
 preserving $\omega_8$. This gives rise to the
 continuous homomorphism
$$\rho'_{\phi}:\Gal(K)\to \OO(\omega_8^{\perp})=W(E_8)$$
such that $\rho_{\phi}$ coincides with the composition of
$\rho'_{\phi}$ and the inclusion map $\OO(\omega_8^{\perp})\subset
\OO(I^{1,8})$. Clearly, $\rho'_{\phi}$ is nothing else but the
defining homomorphism of the Galois module $\omega_8^{\perp}$.

\begin{lem}
\label{omegaC} The Galois modules $J(C)_2$ and
$\omega_8^{\perp}\otimes\Z/2Z$ are isomorphic.
\end{lem}

\begin{proof}
Since both Galois modules $J(C)_2$ and
$\omega_8^{\perp}\otimes\Z/2Z$ are isomorphic to $L/2L$ (Remark
\ref{GCS2} and Sect. \ref{general}), $J(C)_2$ and
$\omega_8^{\perp}\otimes\Z/2Z$ are isomorphic.
\end{proof}


\end{sect}

\begin{thm}
\label{W8P} Let $S$ be a Del Pezzo surface of degree $1$ that is
defined over $K$. Suppose that $\rho'_{\phi}(\Gal(K))=W(E_8)$ or
$W^{+}(E_8)$. Then $\End(J(C))=\Z$.
\end{thm}

\begin{proof}
Recall (Sect. \ref{weyl}) that
$W^{+}(E_8)=\SO(\omega_8^{\perp})=\SO(Q(E_8))$. Replacing (if
necessary) $K$ by its suitable quadratic extension, we may and
will assume that $\rho'_{\phi}(\Gal(K))=W^{+}(E_8)$.
 By Lemma \ref{omegaC},
$J(C)_2$ and $\omega_8^{\perp}\otimes\Z/2Z$ are isomorphic. In
light of Theorem \ref{AVO8}, in order to finish the proof, it
suffices to check that the image of $W^{+}(E_8)$ in $\Aut(
Q(E_8)\otimes\Z/2\Z)$ is isomorphic to $\OO^{+}(8)$. But this
assertion follows easily from
 \cite[Ch. VI, Sect. 4, Ex. 1]{Bourbaki}.
\end{proof}

\begin{thm}
\label{WA8} Let $S$ be a Del Pezzo surface of degree $1$ that is
defined over $K$. Suppose that $\rho'_{\phi}(\Gal(K))$ contains a
subgroup that is a conjugate of $\iota(\A_{8})$ in  $W(E_{8})$.

Then $\End^0(J(C))$ is either $\Q$ or a quadratic field. In
particular, $J(C)$ is absolutely simple.
\end{thm}

\begin{proof}
Replacing $\phi$ by its composition with a suitable element of
$W(E_{8})=\OO(\omega_8^{\perp})$, we may and will assume that
$\rho'_{\phi}(\Gal(K))$ contains $\iota(\A_{8})$. Replacing (if
necessary) $K$ by its suitable finite algebraic extension, we may
and will assume that
$$\rho'_{\phi}(\Gal(K))=\iota(\A_{8}).$$
Clearly, the map
$$\rho_{\A}:\Gal(K) \to \A_{8}, \
\sigma\mapsto \iota^{-1}(\rho'_{\phi}(\sigma))$$ is a continuous
{\sl surjective} group homomorphism and
$$\rho'_{\phi}=\iota \ \rho_{\A}.$$
Let us consider the $8$-element set $\B=\{1,2 \dots, 8\}$ and the
 $8$-dimensional vector space $\F_2^8=\F_2^{\B}$ provided with the
 natural structure of $\A_8$-module. The surjective homomorphism
 $\rho_{\A}:\Gal(K)\twoheadrightarrow \A_8$ provides $\F_2^8$
 with the natural structure of $\Gal(K)$-module; in addition, the
 image of $\Gal(K)$ in $\Aut(\F_2^8)$ coincides with $\A_8$. It
 follows from Lemma \ref{F8S} that the Galois modules
 $\omega_8^{\perp}\otimes\Z/2\Z$ and $\F_2^8$ are isomorphic. On
 the other hand, thanks to Lemma \ref{omegaC},
$J(C)_2$ and $\omega_8^{\perp}\otimes\Z/2Z$ are isomorphic.
 This
implies that the Galois modules $\F_2^{\B}=\F_2^8$ and $J(C)_2$
are isomorphic. Now the assertion follows from Theorem
\ref{mainAV}.
\end{proof}

\begin{sect}
\label{generalB}

Let $B \subset \P^2(\bar{K})$ be an 8-element set of points on the
projective plane that enjoys the following properties:

\begin{itemize}
\item[(i)] The set $B$ is {\sl general position}, i.e., no $3$
points  lie on a line, no $6$ points line on a conic, $B$ does not
lie on a singular cubic in such a way that the singular point
belongs to $B$.

\item[(ii)] The group $\Gal(K)$ permutes transitively points of
$B$ and therefore we get a natural homomorphism $\rho_B:\Gal(K)
\to \Perm(B)\cong \mathrm{S}_8$ from $\Gal(K)$ to the group
$\Perm(B)$ of permutations  of $B$.
\end{itemize}
We write $G_B$ for the image of $\Gal(K)$ in $\Perm(B)$ and
consider the $8$-dimensional $\F_2$-vector space $\F_2^B$ of all
$\F_2$-valued functions on $B$ provided with the natural structure
of Galois module.
Let $K(B)\subset \bar{K}$ be the smallest extension of $K$ over which every point
of $B$ is defined. Clearly,  $K(B)/K$ is a finite Galois extension that corresponds to
the kernel of $\Gal(K)\twoheadrightarrow G_B$ and $\Gal(K(B)/K)$ is canonically isomorphic to $G_B$.


 Let $S_B$ be surface that is obtained by blowing up
points of $B$. Then $S_B$ is defined over $K$. Since $B$ is in
general position, $S_B$ is a Del Pezzo surface of degree $1$
 (\cite[Sect. 3]{Iskovskikh}, \cite[Th. 1 on p. 27]{Demazure}).
 We write $C_B$ for the corresponding branch curve.

\end{sect}

\begin{lem}
\label{order2} The Galois modules $J(C_B)_2$ and $\F_2^B$ are
isomorphic.  In particular, $K(J(C_B)_2)=K(B)$, the groups $\tilde{G}_{2,J(C_B),K}$  and $G_B$ are  isomorphic
 and  the field extension $K(J(C_B)_2)/K$  corresponds to the kernel of
 $$\Gal(K)\twoheadrightarrow G_B\subset \Perm(B).$$
\end{lem}

\begin{proof}
Recall that $J(C_B)_2=\Pic(C_B)_2$.
In light of Remark \ref{GCS2}, it suffices to check that the
Galois modules $\F_2^B$ and $L/2L$ are isomorphic.

We mimick the proof of Lemma \ref{F8S}.
 For each $b\in B$ we write $\ell_b$ for the
class in $\Pic(S_B)$ of the corresponding exceptional curve in
$S$. We write $f_0$ for the class  in $\Pic(S_B)$  of the
preimage of line in $\P^2$. Clearly,
$$\sigma(f_0)=f_0, \ \ell_{\sigma(b)}=\sigma(\ell_{b}) \ \forall b\in B,
\sigma\in\Gal(K).$$ It is known \cite[Sect. 25.1.2]{Manin} that
$f_0$ and $\{\ell_b\}_{b\in B}$ constitute a basis of the free
commutative group $\Pic(S_B)$,
$$K_{S_B}= -3f_0+\sum_{b\in B}\ell_B\in \Pic(S_B),$$
$$(f_0,f_0)=1, (\ell_b,\ell_b)=-1, (\ell_b,f_0)=0 \ \forall b
\in B$$ and $(\ell_{b_1},\ell_{b_2})=0$ if $b_1\ne b_2$.

Recall that $(\K_{S_B},\K_{S_B})=1$. It follows that $\Pic(S_B)$
splits into orthogonal Galois-invariant direct sum
$$\Pic(S_B)=\Z\cdot \K_{S_B}\oplus L.$$
We write $\bar{f}_0$ and $\bar{\ell}_b$ for the images in
$\Pic(S_B)/2\Pic(S_B)$ of $f_0$ and $\ell_b$ respectively.
Clearly, the image of $K_{S_B}$ in $\Pic(S_B)/2\Pic(S_B)$
coincides with
$$\bar{\omega}:=\bar{f}_0+\sum_{b\in B}\bar{\ell}_b;$$
notice that $\bar{f}_0$ and $\{\bar{\ell}_b\}_{b\in B}$ constitute
a basis of the $\F_2$-vector space $\Pic(S_B)/2\Pic(S_B)$. The
unimodular intersection pairing on $\Pic(S_B)$ induces a
non-degenerate pairing
$$\Pic(S_B)/2\Pic(S_B)\times \Pic(S_B)/2\Pic(S_B) \to \F_2.$$
Clearly, the splitting
$$\Pic(S_B)/2\Pic(S_B)=\F_2\cdot \bar{\omega}\oplus L/2L$$
is Galois-invariant and orthogonal. In particular, $L/2L$
coincides with the orthogonal complement of $\bar{\omega}$
 in $\Pic(S_B)/2\Pic(S_B)$. It follows that
$$L/2L=\{c \bar{f}_0+\sum_{b\in B}c_b\bar{\ell}_b\mid c+\sum_{b\in
B}c_b=0\}.$$ Now the map $$c \bar{f}_0+\sum_{b\in
B}c_b\bar{\ell}_b \to \{b \mapsto c_b\}$$ establishes an
isomorphism of Galois modules $L/2L$ and $\F_2^B$.
\end{proof}

Let us consider the jacobian $J(C_B)$ of $C_B$; it is a
four-dimensional abelian variety defined over $K$.

\begin{thm}
\label{main}
Suppose that $G_B=\Perm(B)$ or $\Alt(B)$. Then:
\begin{itemize}
\item[(i)]
 $\End^0(J(C_B))$ is either $\Q$ or a quadratic field.
In particular, $J(C_B)$ is absolutely simple.
\item[(ii)]
Let $T \subset \P^2(\bar{K})$ be another Galois-invariant set of $8$ points in general position.
Assume that either $G_T$ is a solvable group or $K(S)$ and $K(T)$ are linearly disjoint over $K$.
Then the abelian varieties $J(C_B)$ and $J(C_T)$ are not isomorphic over $\bar{K}$.
\end{itemize}
\end{thm}

\begin{proof}
 Replacing (if necessary) $K$ by
suitable quadratic extension, we may and will assume that
$G_B=\Alt(B)$.
Now, combining Theorem \ref{mainAV} and Lemma
\ref{order2}, we obtain the assertion (i). Applying Theorem \ref{nonisom} to $X=J(C_B)$ and
$Y-J(C_T)$, we obtain the assertion (ii).
\end{proof}

 The next statement explains how to construct points in
general position.

\begin{prop}\label{General}
Suppose that $E\subset \P^2$ is an absolutely irreducible cubic
curve that is defined over $K$. Suppose that $B\subset E(K_a)$ is
a  $8$-element set that is a $\Gal(K)$-orbit. Let us assume that
 the image $G_B$ of $\Gal(K)$ in the group $\Perm(B)$
of all permutations of $B$ coincides either with $\Perm(B)$ or
with  $\Alt(B)$. Then $B$ is in general position.
\end{prop}

\begin{proof}
Notice that $\Gal(K)$ acts $3$-transitively on $B$.

Step 1. Suppose that $D$ is a line in $\P^2$  that contains three
points of $B$ say,
$$\{P_1,P_2,P_3\}\subset \{P_1,P_2,P_3, P_4,P_5,P_6,P_7,P_8\}=B.$$
Clearly, $D\bigcap E=\{P_1,P_2,P_3\}$. There exists $\sigma
\in\Gal(K)$ such that $\sigma(\{P_1,P_2,P_3\})=\{P_1,P_2,P_4\}$.
It follows that the line $\sigma(D)$ contains $\{P_1,P_2,P_4\}$
and therefore $\sigma(D)\bigcap E=\{P_1,P_2,P_4\}$. In particular,
$\sigma(D)\ne D$. However, the distinct lines $D$ and $\sigma(D)$
meet each other at {\sl two} distinct points $P_1$ and $P_2$.
Contradiction.

Step 2. Suppose that $Y$ is a conic in $\P^2$ such that $Y$
contains six points of $B$ say, $\{P_1,P_2,P_3,
P_4,P_5,P_6\}=B\setminus\{P_7, P_8\}$. Clearly, $Y\bigcap
E=B\setminus\{P_7,P_8\}$. If $Y$ is reducible, i.e., is a union of
two lines $D_1$ and $D_2$ then either $D_1$ or $D_2$ contains (at
least) three points of $B$, which is not the case, thanks to Step
1. Therefore $Y$ is {\sl irreducible}.

There exists $\sigma \in\Gal(K)$ such that $\sigma(P_1)=P_7,
\sigma(P_8)=P_8$. Then $\sigma(P_7)=P_i$ for some positive integer
$i\le 6$. This implies that
$\sigma(B\setminus\{P_7,P_8\})=B\setminus\{P_i,P_8\}$ and
 the irreducible conic $\sigma(Y)$ contains
$B\setminus\{P_i,P_8\}$. Clearly, $\sigma(Y)\bigcap
E=B\setminus\{P_i,P_8\}$ contains $P_7$. In particular,
$\sigma(Y)\ne Y$. However, both conics contain the $5$-element set
$B\setminus\{P_i,P_7,P_8\}$. Contradiction.

Step 3. Suppose that $Z$ is a cubic in $\P^2$ such that $B\subset
Z$ and say, $P_1\in B$ is a singular point of $Z$. If $Z$ is
reducible then either there is a line with $3$ points of $B$ or a
conic with $6$ points of $B$. So, $Z$ is irreducible and therefore
$P_1$ is the only singular point of $Z$. Clearly, for each
$\sigma\in \Gal(K)$ the cubic $\sigma(Z)$ also contains $B$ and
$\sigma(P_1)$ is the only singular point of $\sigma(Z)$. Pick
$\sigma$ with $\sigma(P_1)=P_2\in B$. Then $\sigma(Z)\ne Z$. The
cubics $Z$ and $\sigma(Z)$ meet at all $8$ points of $B$. In
addition, the local intersection index at  singular $P_1$ and
$\sigma(P_1)=P_2$ is, at least, $2$. This implies that the
intersection index of $Z$ and $\sigma(Z)$ is, at least $10$, which
is not true, since the index is $9$.
\end{proof}

\begin{proof}[Proof of Theorem \ref{mainex}]
Let $f(t)\in K[t]$ be an irreducible  degree $8$ polynomial, whose Galois
group $\Gal(f)$ is either $\ST_8$ or $\A_8$. Let $\RR\subset K_a$ be the set
of roots of $f$. Then the $8$-element set
$$B(f)=\{(\alpha^3:\alpha:1)\in \P^2(K_a)\mid \alpha\in \RR\}$$
lies on the absolutely irreducible $K$-cubic $xz^2-y^3=0$  and $G_{B(f)}\cong
\ST_8$ or $\A_8$ respectively. It follows from Proposition \ref{General} that
$B(f)$ is in general
 position, which proves (i).  Clearly, $K(B(f))$ coincides with the {\sl splitting field}
 $K(\RR)$ of $f$. Therefore
$$G_{B(f)}=\Gal(K(B(f))/K)=\Gal(K(\RR)/K)= \Gal(f),$$
 which implies that $G_{B(f)}$ is either
  $\Perm(B)$ or $\Alt(B)$. Now the assertions (ii) and (iii) follow from the first and second assertions
 of Theorem \ref{main} respectively.


\end{proof}

\begin{ex}
Let $L=\C(y_1, \dots ,y_8)$ be the field of rational
functions in $8$ independent variables over $\C$.  There is the natural  action of $\ST_8$ on $L$ by $\C$-linear field automorphisms such that each permutation $s$ sends every $y_i$ to $y_{s(i)}$.
Let $K$ be the subfield of $\ST_8$-invariants of $L$. Clearly, $\C\subset L$, the field extension $L/K$ is a finite Galois extension,
$\Gal(L/K)=\ST_8$ and $\bar{L}=\bar{K}$. Let us choose $c\in K$ (e.g., $c=0$ or $c=(\sum_{i=1}^8 y_i)/8$) and consider the $8$-element set
$$B=\{(y_i-c)^3:(y_i-c):1)\mid 1\le i\le 8\}\subset \P^2(L)\subset  \P^2(\bar{L})=  \P^2(\bar{K}).$$
Clearly, $B$ is $\Gal(K)$-invariant and $G_B=\Perm(B)$. It follows that $B$ is in
{\sl in general position}.

Applying
Theorem \ref{main}, we conclude that $\End^0(J(C_{B}))$ is either $\Q$ or a quadratic field; in particular,
 $J(C_{B})$ is an (absolutely) simple abelian variety.
\end{ex}

\section{Explicit formulas}
\label{explicit}
Let $h(t)\in K[t]$ be an irreducible degree $8$ polynomial, whose Galois group is either $\ST_8$ or $\A_8$. In order to  simplify slightly our computations, we assume that $h(t)$ is {\sl monic} and the sum of its roots is zero, i.e.,
$$h(t)=t^8+\sum_{i=0}^6 h_i t^i; \ h_i \in K.$$
The irreducibility of $h(t)$ implies that $h_0\ne 0$.
Let
$K[t]_7$ be the $8$-dimensional subspace in
$K[t]$ that consists of all polynomials, whose degree does not exceed $7$. We write
$$D_h:K[t] \twoheadrightarrow K[t]_7$$
for the surjective $K$-linear map that sends any polynomial into its remainder with respect to division by $h(t)$. By definition,
$f-D_h(f)$ is divisible by $h$ for all $f\in K[t]$.

Let $K[x]$ and $K[x,y]$ be the ring of polynomials in independent variables $x$ and $x,y$ respectively. We write
$$A:K[t] \to K[x]\oplus y\cdot K[x]\oplus y^2\cdot K[x]\subset K[x,y]$$
for the $K$-linear map that sends $t^{3i}$ to $x^i$, $t^{3i+1}$ to $x^i y$ and $t^{3i+2}$ to
$x^i y^2$ respectively (for each nonnegative integer $i$).

Clearly, if $g(t)\in K[t]$ and $G(x,y)=A(g)$ then $g(t)=G(t^3,t)$ and $G(x,y)-g(y)$ is divisible by $x-y^3$ in $K[x,y]$. In addition, $\deg(G)\ge \deg(g)/3$. On the other hand,
 if $\deg(g)\le 3d$
for some positive integer $d$ and $g$ does {\sl not} have a term of degree $3d-1$  then $\deg(G)\le d$. For example, if $\deg(g)=9$ and $g$ does not contain term $t^8$ then $G$ contains term $x ^3$  and $\deg(G)= 3$. Another examples: a)if $\deg(g)\le 7$ then $\deg(G) \le 3$;
b)if $\deg(g)=16$ then $G$ contains term $x^5 y$ and $\deg(G)=6$.

Let us put
$$B=B(h).$$
Our first goal is to describe explicitly the (two-dimensional) $\bar{K}$-vector space $H$ of all cubic forms (in homogeneous coordinates $(x:y:z)$ and with coefficients in $\bar{K}$) that do vanish at all points of $B$ and find the ``nineth point". Clearly, one of those forms is $u:=xz^2-y^3$.
 Another  cubic form $v(x,y,z)$ is defined by
 $v(x,y,1)=A(th(t))$.
Since $v(x,y,1)-yh(y)=A(th(t))-yh(y)$ is divisible by $x-y^3$, the form $v$
does vanish at all points of $B$. Since $th(t)$ has degree $9$, the polynomial $ A(th(t))$ contains the term $x^3$ and therefore has $x$-degree $3$. It follows that $v$ also has $x$-degree $3$. Since $u$ and $v$ have different $x$-degrees, they are not proportional one to another and therefore
$$H=\bar{K}\cdot u + \bar{K}\cdot v$$ is a two-dimensional space of cubic forms.
Clearly, both $u$ and $v$ do vanish at  $(0:0:1)$: this is the nineth point \cite[Ch. 5, Sect. 4, Cor. 4.5]{H} and the set of common zeros of $u$ and $v$ coincides with $B\cup \{(0:0:1)\}$.

Recall that $S_{B}$ is obtained from $\P^2$ by blowing up points of $B$. We write
$$g_B: S_{B} \to \P^2$$
for the corresponding regular birational map. For each $b\in B$  its preimage $E_b:=g_B^{-1}(b)$ is a smoth projective rational curve with self-intersection index $-1$. By definition, $g_B$ establishes a biregular isomorphism between $S_B\setminus\bigcup_{b\in B}E_b$ and $\P^2\setminus B$. Clearly, the line $L:z=0$ does not meet $B$. Further, we view $L$ as a divisor on $\P^2$. It is known \cite[Sect. 25.1 and 25.1.2 on pp. 126--127]{Manin} that
$$\K_{S_B}:=-3 g_B^{*}(L)+\sum_{b\in B}E_b= -g_B^{*}(3L)+\sum_{b\in B}E_b$$
is a {\sl canonical} divisor on $S_B$. Clearly, for each form $q\in H$ the rational function $q/z^3$ on $\P^2$ satisfies $\div(q/z^3)+3L\ge 0$, i.e., $q/z^3\in \Gamma(\P^2,3L)$. Also $q/z^3$ is defined and vanishes at every point in $B$. It follows that $q/z^3$ (viewed as a rational function on $S_B$) lies in
$\Gamma(S_B, 3 g_B^{*}(L)-\sum_{b\in B}E_b)=\Gamma(S_B,-\K_{S_B})$. Since the latter space is two-dimensional,
$$\Gamma(S_B,-\K_{S_B})=\bar{K}\cdot\frac{u}{z^3}\oplus \bar{K}\cdot\frac{v}{z^3}.$$
This gives us a rational anticanonical map
$$S_B \stackrel{g_B}{\longrightarrow}\P^2 \stackrel{(u:v)}{\longrightarrow} \P^1,$$
which is regular outside the (preimage of the) base point $(0:0:1)$.

Now let us consider the space $H_2$ of degree $6$ forms that vanish with all its  partial derivatives of first order at every point of $B$. For example,
$$u^2, uv, v^2 \in H_2;$$
they are linearly independent over $\bar{K}$, because they have (distinct) $x$-degrees $2,4,6$ respectively. (Notice that all of them do vanish at $(0:0:1)$.) One may easily check (using, for instance,  \cite[Ch. IV, Sect. 3.1]{Sha}) that if $F\in H_2$ then $F/z^6$ (viewed as a rational function on $S_B$) lies in
$\Gamma(S_B, 6 g_B^{*}(L)-\sum_{b\in B}2E_b)=\Gamma(S_B,-2\K_{S_B})$. Recall \cite[p. 68]{Demazure} that the space $\Gamma(S_B,-2\K_{S_B})$ is four-dimensional. This implies that if a form $w\in H_2$ does not vanish at $(0:0:1)$ then it is {\sl not} a linear combination of $u^2, uv, v^2$ and therefore $\{u^2, uv, v^2, w\}$ is a basis of $H_2$ and
$$\Gamma(S_B,-2\K_{S_B})=\bar{K}\cdot \frac{u^2}{z^6}\oplus \bar{K}\cdot
\frac{uv}{z^6}\oplus\bar{K}\cdot \frac{v^2}{z^6}\oplus \bar{K}
\cdot \frac{w}{z^6}.$$
This gives us a regular doubly-anticanonical map
$$\pi_B:S_B \stackrel{g_B}{\longrightarrow}\P^2 \stackrel{(u^2:uv:v^2:w)}{\longrightarrow} \P^3.$$
Recall that $C_B$ is the branch curve of $\pi_B$. This allows us to describe explicitly the
plane projective curve $g_B(C_B)$, which is birationally isomorphic to $C_B$. Indeed,
$g_B(C_B)\setminus B$ coincides with the set of points of $\P^2\setminus B$ where the the tangent map to
$$(x:y:z) \mapsto (u(x,y,z)^2:u(x,y,z) v(x,y,z):v(x,y,z)^2:w(x,y,z))$$
is {\sl not} injective. Since $B\bigcup \{(0:0:1)\}$ is the set of common zeros of $F_1$ and $F_2$, the curve $g_B(C_B)\setminus B\setminus \{(0:0:1)\}$ coincides with the set of points of $\P^2\setminus B\setminus \{(0:0:1)\}$ where the form
$$Q(x,y,z):=
\begin{vmatrix}
u_x & u_y & u_z\\
v_x & v_y & v_z\\
w_x & w_y & w_z
\end{vmatrix}
$$
does vanish. Since $C_B$ is an irreducible projective curve, $g_B(C_B)$ is also projective ireducible and coincides with the set of zeros of $Q$, i.e.,
$$g_B(C_B)=\{(x:y:z)\in \P^2\mid Q(x,y,z)=0\}\subset\P^2.$$
Clearly, $Q$ has coefficients in $K$ and degree $9$. Since $g_B(C_B)$ is irreducible and $C_B$ has genus $4$,
the curve $C_B$ is singular and $Q(x,y,z)$ is an irreducible polynomial. Indeed, if $Q$ is reducible then the irreducibility of  $g_B(C_B)$ implies that $Q$ is a power of an irreducible polynomial, i.e., $Q$ is either a $9$th power of a linear form or a cube of a cubic form. In the former case, $g_B(C_B)$ has genus $0$ while in the latter one its genus is either $0$ or $1$.
However, $g_B(C_B)$ is birationally isomorphic to the genus $4$ curve $C_B$, which rules out both possibilities for the reducibility of $Q$. Since $\deg(Q)=9$ and $4 \ne (9-1)(9-2)/2$, the curve $g_B(C_B)$ is singular. Clearly, all its singular points must lie in $B$. Since $Q$ has coefficients in $K$ and $B$ constitutes a Galois orbit,  all points of $B$ are singular
and have the same multiplicity. A well-known formula for the genus of (the normalization of) a plane curve \cite[Ch. IV, Sect. 4.1]{Sha} implies that all points of $B$ have multiplicity $3$.

The rest of this Section is devoted to an explicit construction of $w$. Let us consider the polynomial
$$h(t)^2=t^{16}+\sum_{i=0}^{14}c_i t^i, \ c_i\in K, c_0=h_0^2\ne 0.$$
Clearly, all the coefficients $c_i$ can be expressed explicitly in terms of the coefficients of $h(t)$. Let us put $F(x,y):=A(h(t)^2)$. Clearly, $\deg(F)=6$ and $F(x,y)-h(y)^2$ is divisible by $x-y^3$. In particular, $F$ does vanish at all points $(\alpha^3,\alpha)$ where $\alpha$ is a root of $h(t)$. It also follows that $3y^2 F_x+F_y$ does vanish at $(\alpha^3,\alpha)$. (On the other hand, $F$ does {\sl not} vanish at $(0,0)$.) If $G(x,y)$ is any polynomial then clearly, both $H(x,y)=F(x,y)-G(x,y)(x-y^3)$ and $3y^2 H_x+H_y$ do vanish
at $(\alpha^3,\alpha)$ while $H$ does {\sl not} vanish at $(0,0)$. I want to find such $G$ that $\deg(G)\le 3$ and $H_x$ does vanish at all $(\alpha^3,\alpha)$. Since $3y^2 H_x+H_y$ do vanish
at $(\alpha^3,\alpha)$, we conclude that both $H_x$ and $H_y$ do vanish
at each $(\alpha^3,\alpha)$. After that, we define
the degree $6$ form $w$ by $w(x,y,1):=H(x,y)$, i.e.
$$w(x,y,z)=z^6 H(x/z,y/z).$$
Clearly, $w$ vanishes at all points of $B=B(h)$ with its partial derivatives of first order with respect to $x$ and to $y$.
By Euler's theorem,
$$x w_x+y w_y+ z w_z=6 w.$$
It follows that $w_z$ also does vanish at all points of $B$. However, $$w(0,0,1)=H(0,0)\ne 0.$$
So, such $w$ enjoys the desired properties and now
our task boils down to finding such $H$ (i.e., finding such $G$).

Now let us put
$$p(t):=D_h(F_x(t^3,t))\in K[t]_7\subset K[t], \ G(x,y): = A(p(t))=A (D_h(F_x(t^3,t))).$$
I claim that $\deg(G)\le 3$ and $G(\alpha^3,\alpha)=F_x(\alpha^3,\alpha)$.
Indeed, $\deg(p(t))\le 7$ and therefore $G(x,y)= A (p(t))$ has degree, at most $3$. Since $F_x(t^3,t)-D_h(F_x(t^3,t))$ is divisible by $h(t)$ and $h(\alpha)=0$,
the values of $F_x(t^3,t)$ and $p(t)$ at $t=\alpha$ do coincide, i.e., $F_x(\alpha^3,\alpha)=p(\alpha)$. On the other hand,
$G(x,y)-p(y)$ is divisible by $x-y^3$ and therefore $G(\alpha^3,\alpha)=p(\alpha)$.
We have
$$G(\alpha^3,\alpha)=p(\alpha)=F_x(\alpha^3,\alpha).$$
If we put (as above) $H:=F-(x-y^3)G$ then
$$H_x(\alpha^3,\alpha)=F_x(\alpha^3,\alpha) -\{G_x(\alpha^3,\alpha)\cdot (\alpha^3-\alpha^3)+G(\alpha^3,\alpha)\cdot 1\}=F_x(\alpha^3,\alpha)-G(\alpha^3,\alpha)=0.$$
We are done.

\end{document}